\documentclass[11pt, a4paper]{article}
\usepackage{amsmath,amssymb,latexsym,graphics,epsfig}
\usepackage{hyperref}
\usepackage{color}
\usepackage{amsthm}
\usepackage{cases}

\newtheorem{theorem}{Theorem}[section]

\def\prova{{\boldmath  $Proof.$}\hskip 0.3truecm}
\def\final{\mbox{ \quad $\Box$}}
\def\DS{\mbox{{\tiny $\!D\!S$}}}
\def\MH{\mbox{{\tiny $\!M\!H$}}}
\def\NA{\mbox{{\tiny $\!N\!A$}}}

\begin{document}

\title{The $(\Delta,D)$ and $(\Delta,N)$ Problems for \\
New Amsterdam and Manhattan Digraphs}
\author{C. Dalf\'o and M.A. Fiol
\\ \\
{\small Universitat Polit\`ecnica de Catalunya, BarcelonaTech} \\
{\small Dept. de Matem\`atica Aplicada IV} \\
{\small Barcelona, Catalonia} \\
{\small (e-mails: {\tt \{cdalfo,fiol\}@ma4.upc.edu})}}

\date{}

\maketitle

\begin{flushright}
To Joan Gimbert, in memoriam.
\end{flushright}

\begin{abstract}
We give a quasi-complete solution of the $(\Delta,N)$ problem
for two well-known families of digraphs used as good models for large interconnection networks. In our study we also relate both families, the New Amsterdam and Manhattan digraphs, with the double-step graphs (or circulant graphs with degree two).
\end{abstract}

{\small \textbf{Keywords:} $(\Delta,D)$ problem, $(\Delta,N)$ problem, New Amsterdam digraph, Manhattan digraph}

\indent {\small {\bf 2000 MSC:} 05C20,05C12}

\section{Introduction}
In this paper we concentrate on two families of digraphs which are Cayley
digraphs of the plane crystallographic groups. Namely, we consider
the $(\Delta,D)$ and $(\Delta,N)$ problems for the so-called New Amsterdam and
Manhattan digraphs. The first problem consists of maximizing the
number $N$ of vertices giving the (maximum degree) $\Delta$ and the
diameter $D$, whereas the second one (somehow dual of the first)
consists of minimizing the diameter for a fixed degree and number of
vertices. Our study is based on known results about the same
problems of another family, the double-step graphs (also called
circulants or Cayley graphs of Abelian groups). In fact, the
$(\Delta,D)$ problem was already solved by Morillo, Fiol and
F\`abrega in~\cite{mff85} for the New Amsterdam digraphs with odd
diameter, and also for the Manhattan digraphs with even diameter.
Although a solution for the other diameters was also claimed in the
same paper, the digraphs proposed were not vertex transitive and, as
a consequence, the eccentricity from the odd vertices was not the
correct one. Here we show that, for such values of the diameter, the
number of vertices is much smaller than the theoretical upper
(Moore-like) bounds. For a comprehensive survey on Moore graphs and the $(\Delta,D)$ problem, see Miller and Sir\'an
\cite{ms05}.

\section{Double-step graphs, New Amsterdam and Manhattan digraphs}
In this section we define the different families of digraphs considered and recall the
corresponding theoretical (Moore-like) upper bounds for their number of vertices.

\subsection{Double-step graphs}
A {\em double-step graph} $G(N;\pm a,\pm b)$ has set of vertices $\mathbb{Z}_N$ 
(the integers modulo $N$) and every vertex $i$ is adjacent to the vertices $i\pm a$ 
and $i\pm b$ (arithmetic is always $\textrm{mod } N$), for some different integers $a,b$ 
called {\em steps} such that $\gcd(N,a,b)=1$.
For more details, see Yebra, Fiol, Morillo and Alegre~\cite{yfma85}.
Then, as is readily seen in that paper, the maximum number $N_{\DS}$ of vertices of a double-step graph
with diameter $k$ is upper bounded  by the Moore-like bound
\begin{equation}
\label{Moore-DS}
N_{\DS}\leq M_{\DS}(2,k)=1+\sum_{n=1}^k 4n=2k^2+2k+1.
\end{equation}

\subsection{New Amsterdam digraphs}
\begin{figure}[h]
\begin{center}
\includegraphics[width=11cm]{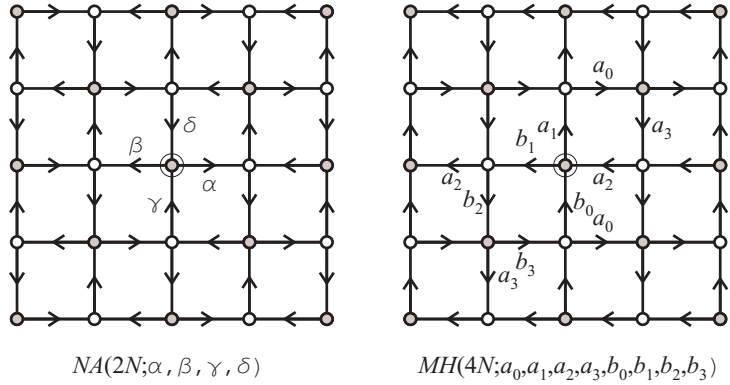}
\vskip-9cm
\caption{A New Amsterdam digraph (left) and a Manhattan digraph (right). In both cases, the even vertices are grey and the odd ones are white.}
\label{fig-new-amsterdam-manhattan}
\end{center}
\end{figure}

Let $N$ be an even integer and let $\alpha,\beta,\gamma,\delta$ be some odd integers ($\alpha\neq\beta$) as before
called {\em steps} satisfying
\begin{equation}\label{cond-stepsNA}
\alpha+\beta+\gamma+\delta \equiv 0\ (\textrm{mod } N).
\end{equation}
Then, a {\em New Amsterdam digraph} $N\!A(N;\alpha,\beta,\gamma,\delta)$ is a bipartite digraph with
set of vertices $V=\mathbb{Z}_N$, $V=V_0\cup V_1$,
$V_0=\{0,2,\ldots,N-2\}$, $V_1=\{1,3,\ldots,N-1\}$, and where each vertex $i\in V_0$ is adjacent to the vertices $i+\alpha,i+\beta\in V_1$, and every vertex $j\in V_1$ is adjacent to the vertices $j+\gamma,j+\delta\in V_0$.
See Figure~\ref{fig-new-amsterdam-manhattan} (left) for the plane local pattern followed by the vertices.

Since the digraph is regular and bipartite, if it has diameter $k$, its maximum number $N_k$ of vertices is twice the number of vertices in $V_0$ when $k$ is odd, or in $V_1$ if $k$ is even, at distance at most $k-1$ from vertex $0$. This leads to the following Moore-like bounds (see Morillo, Fiol, and F\`{a}brega \cite{mff85}):
\begin{numcases}{N_{\NA}\leq M_{\NA}(2,k)=}
2\left(1+\sum_{m=1}^n 4m\right)=k^2+1 & \hskip-.25cm for odd $k=2n+1$,\label{Moore-NA-odd}\\
2\sum_{m=1}^n (4m-2)=k^2              & \hskip-.25cm for even $k=2n$.\label{Moore-NA-even}
\end{numcases}

\subsection{Manhattan digraphs}
Manhattan digraphs were introduced independently (in slightly different forms and contexts) by Morillo, Fiol, and F\`abrega \cite{mff85} and Maxemchuk \cite{m87}. They are called in this way because locally resemble the topology of the avenues and streets of Manhattan (or {\em l'Eixample} in downtown Barcelona); see Figure~\ref{fig-new-amsterdam-manhattan} (right).
More precisely, given an integer $N$ multiple of $4$, the {\em Manhattan digraph} $M\!H=M\!H(N;a_0,b_0,a_1,b_1,a_2,b_2,a_3,b_3)$ has set of vertices
$V=\bigcup_{j=0}^3 V_j$, where $V_j=\{i\in\mathbb{Z}_N : i\equiv j\ (\textrm{mod } 4)\}$. Moreover, each vertex $i\in V_j$ is adjacent to the vertices $i+a_j$, $i+b_j$, where
the steps $a_0,b_0,\ldots,a_3,b_3$ satisfy $a_j\equiv 3$ $(\textrm{mod } 4)$, $b_j\equiv 1$ $(\textrm{mod } 4)$, for $j=0,1,2,3$, and
\begin{equation}\label{cond-steps-MH}
 a_0+a_2 \equiv -(a_1+a_3)\equiv b_0+b_2 \equiv -(b_1+b_3)\ (\textrm{mod } N).
\end{equation}

The number of vertices for such a digraphs
are upper bounded by the following Moore-like bounds (see again Morillo, Fiol, and F\`{a}brega \cite{mff85}):
\begin{numcases}{N_{\MH}\leq M_{\MH}(2,k)=}
2(k-1)^2 & for odd $k=2n+1$, \label{Moore-MH-odd}\\
2[(k-1)^2+1] & for even $k=2n$.\label{Moore-MH-even}
\end{numcases}

\section{New Amsterdam and Manhattan digraphs associated to double-step graphs}
In this section we show that every double-step graph has associated both a New Amsterdam digraph and a Manhattan digraph. In each step of this process, while the number of vertices is doubled, the diameter only increases (at most) by one. Moreover, the relationships between the steps characterizing the digraphs are easily devised by superimposing the corresponding plane patterns.

\subsection{From a double-step graph to a New Amsterdam digraph}
\begin{figure}[h]
\begin{center}
\includegraphics[width=12cm]{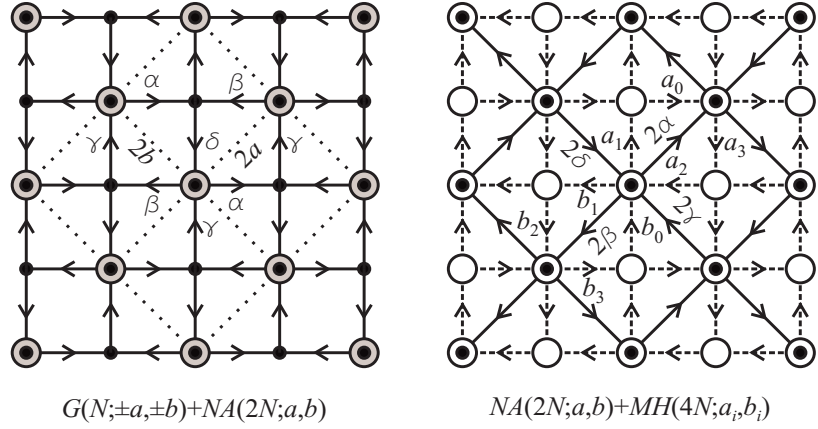}
\vskip-9.5cm
\caption{Left: A New Amsterdam digraph $N\!A(2N;a,b)$ superimposed on a double-step graph $G(N;\pm a,\pm b)$. Right: A Manhattan digraph $M\!H(4N;a_i,b_i), \ 0\leq i\leq3$, superimposed on a New Amsterdam digraph $N\!A(2N;\alpha,\beta,\gamma,\delta)$. The double-step digraph has grey vertices and its edges are dotted lines, the New Amsterdam digraph has black vertices and continuous arcs, and the Manhattan digraph has white vertices and discontinuous arcs.}
\label{fig-graf-doble-pas+NA+M}
\end{center}
\end{figure}

Given a double-step graph $G(N;\pm a,\pm b)$, we can obtain a New Amsterdam digraph $N\!A(N_{\NA};\alpha,\beta,\gamma,\delta)$, where $N_{\NA}=2N$ and the steps $\alpha,\beta,\gamma,\delta$ satisfy the following conditions:
\begin{itemize}
\item[$(i)$]
$\alpha,\beta,\gamma,\delta$ are odd integers,
\item[$(ii)$]
$\alpha+\gamma\equiv -\beta-\delta \equiv 2a\ (\textrm{mod } N_{\NA})$,
\item[$(iii)$]
$\beta+\gamma\equiv -\alpha-\delta \equiv 2b\ (\textrm{mod } N_{\NA})$,
\end{itemize}
where $(ii),(iii)$ are made clear in Figure~\ref{fig-graf-doble-pas+NA+M}.
There are different possible solutions of these equations. For instance, we can take
\begin{equation}\label{ds->na}
\alpha=-1,\qquad \beta=2(b-a)-1,\qquad \gamma=2a+1,\qquad \delta=-2b+1,
\end{equation}
and we are lead to the following result.

\begin{theorem}
\label{basic-theo}
For any double-step graph $G=G(N;\pm a,\pm b)$ with diameter $D=k$, there exists a New Amsterdam
digraph $N\!A=N\!A(N_{\NA};\alpha,\beta,\gamma,\delta)$ on $N_{\NA}=2N$ vertices, steps
$\alpha=-1,\beta=2(b-a)-1, \gamma=2a+1, \delta=-2b+1$, and diameter $D_{\NA}$ satisfying
\begin{equation}\label{diameterNA}
2k\le D_{\NA}\le 2k+1.
\end{equation}
\end{theorem}
\prova
Let us first prove that there is a path of length at most $2k+1$ between any two given vertices
$i,j$ of a New Amsterdam digraph $N\!A$. By symmetry, it is enough to check the cases $i\in \{0,1\}$. Then we must
distinguish the cases when $j$ is even or odd.
\begin{itemize}
\item[$(a)$] $0\rightarrow2j$:
Since $\gcd(a,b,N)=1$ and $G$ has diameter $k$, for any $j\in \mathbb{Z}_N$
there exist integers $m,n$ such that $|m|+|n|\le k$ and
\begin{equation}
\label{ma+nb}
ma+nb \equiv j\ (\textrm{mod } N).
\end{equation}
So, multiplying by 2 and using Eqs. $(ii)$ and $(iii)$, we have that
\begin{equation}\label{(0,2j)}
m\alpha+m\gamma+n\beta+n\gamma \equiv 2j\ (\textrm{mod } N_{\NA})
\end{equation}
and, hence, $\textrm{dist}_{\NA}(0,2j)\le 2|m|+2|n|\le 2k$.

\item[$(b)$] $0\rightarrow2j-1$:
Since in at most $2k$ steps we reach every even vertex from $0$, by
using an additional step $\alpha=-1$, Eq. (\ref{(0,2j)}) yields
\begin{equation}\label{(0,2j-1)}
(m+1)\alpha+m\gamma+n\beta+n\gamma \equiv 2j-1\ (\textrm{mod } N_{\NA})
\end{equation}
and we reach every odd vertex in at most $2k+1$ steps,
$\textrm{dist}_{\NA}(0,2j-1)\le 2|m|+1+2|n|\le 2k+1$.

\item[$(c)$] $1\rightarrow2j+1$:
By using again Eq. (\ref{(0,2j)}), we have
\begin{equation}\label{(1,2j+1)}
1+m\alpha+m\gamma+n\beta+n\gamma \equiv 2j+1\ (\textrm{mod } N_{\NA})
\end{equation}
and, hence, $\textrm{dist}_{\NA}(1,2j+1)\le 2k$.

\item[$(d)$] $1\rightarrow2j$:
Let us consider the value $2j-\gamma=2(j-a)-1$. By the same reasons
as before, there exist integers $m',n'$ such that $|m'|+|n'|\le k$
and
\begin{equation}
\label{m'a+n'b}
m'a+n'b \equiv j-a-1\ (\textrm{mod } N).
\end{equation}
and, hence,  using again Eqs. $(ii)$ and $(iii)$,
$$
1+m'\alpha+m'\gamma+n'\beta+(n'+1)\gamma\equiv 2j\ (\textrm{mod } N_{\NA})
$$
and, so, $\textrm{dist}_{\NA}(1,2j)\le 2k+1$.
\end{itemize}
Finally to prove that $D_{\NA}\ge 2k$ it suffices to take the vertex
$i$ such that $\textrm{dist}(0,i)=k$. \final

\subsection{From a New Amsterdam digraph to a Manhattan digraph}
Similarly, from a New Amsterdam digraph $N\!A(N_{\NA};\alpha,\beta,\gamma,\delta)$ we can obtain the Manhattan digraph $M\!H(N_{\MH};a_0,b_0,a_1,b_1,a_2,b_2,a_3,b_3)$, with $N_{\MH}=2N_{\NA}$ and  steps $a_0,b_0,\ldots,a_3,$ $b_3$ satisfying:
\begin{itemize}
\item[$(i)$]
$a_i,b_i$ are odd integers for $i=0,1,2,3$.
\item[$(ii)$]
$a_0+a_2 \equiv -(a_1+a_3) \equiv b_0+b_2\equiv -(b_1+b_3)\ (\textrm{mod } N_{\MH})$.
\item[$(iii)$]
$a_0+a_1 \equiv  2\alpha \ (\textrm{mod } N_{\MH})$;\qquad
$b_1+b_2 \equiv  2\beta\ (\textrm{mod } N_{\MH})$; \newline
$b_3-a_1 \equiv  2\delta\ (\textrm{mod } N_{\MH})$; \qquad
$b_0-a_0  \equiv  2\gamma\ (\textrm{mod } N_{\MH})$.
\end{itemize}
For understanding where Eqs. $(ii)$ and $(iii)$ come from, see
Figure~\ref{fig-graf-doble-pas+NA+M} (right).
Again, there are different possible solutions to these equations. For instance, we can take:
\begin{equation}
\begin{array}{llll}
a_0=1, & a_1=2\alpha-1, & a_2=1, & a_3=-2\alpha-1, \\
b_0=2\gamma+1, & b_1=2\beta+2\gamma-1, & b_2=-2\gamma+1, & b_3=-2\beta-2\gamma-1.\label{na->mh2}
\end{array}
\end{equation}

\subsection{From a double-step graph to a Manhattan digraph}
As a consequence of the above results (\ref{ds->na}), and (\ref{na->mh2}), the following equalities give the steps of a Manhattan digraph from the steps of its corresponding double-step graph:
\begin{equation}
\begin{array}{llll}
a_0=1, & a_1=-3, & a_2=1, & a_3=1, \\
b_0=4a+3, & b_1=4b-1, & b_2=-4a-1, & b_3=-4b-1.\label{ds->mh2}
\end{array}
\end{equation}
As a consequence, we get the following theorem.
\begin{theorem}
For any double-step graph $G=G(N;\pm a,\pm b)$ with diameter $D=k$, there exists a Manhattan
digraph $M\!H=M\!H(N_{\MH};a_0,b_0,\ldots,a_3,b_3)$, where $N_{\MH}=4N$ vertices, steps given
by $(\ref{ds->mh2})$, and diameter $D_{\MH}$ satisfying
\begin{equation}\label{diameterNA}
2k+1\le D_{\MH}\le 2k+2.
\end{equation}
\end{theorem}
\prova
As was shown by Dalf\'o, Comellas, and Fiol \cite{cdf08}, every Manhattan digraph can be seen as the line digraph of a New Amsterdam digraph.
Then, the result is a direct consequence of Theorem \ref{basic-theo} and the properties of line digraphs. Namely, if $G$ is a $\delta$-regular digraph (different from a directed cycle) with $N$ vertices and diameter $D$, then its line digraph $L(G)$ has $\delta N$ vertices and diameter $D+1$.
See more details about the line digraph technique in Fiol, Yebra and Alegre~\cite{FiYeAl84}.
\final

\setcounter{section}{4}
\setcounter{equation}{0}
\setcounter{theorem}{0}
\begin{center}
{\bf 4. Dense New Amsterdam and Manhattan digraphs}
\label{four}
\end{center}

Now we are ready to give results about the $(\Delta,N)$ problem for both families considered. We begin recalling some basic known results concerning  double-step graphs.

\subsection{A basic pair of steps}
As it was shown by Bermond, Iliades and Peyrat~\cite{bip85} (see also  Beivide, Herrada, Balc\'azar and Arruabarrena~\cite{bhba91}), the  pair of steps
$$
a=k,\qquad b=k+1,
$$
solves both the
$(\Delta,D)$ and $(\Delta,N)$ problems for the family of double-step graphs.
Indeed, we have the following result.
\begin{theorem}[\cite{bip85,bhba91}]
\label{Bermond-theo}
For any number $N$ of vertices satisfying
$$
M_{\DS}(2,k-1)=2(k-1)^2+2(k-1)+1 < N \le M_{\DS}(2,k)=2k^2+2k+1
$$
the double-step graph $G(N;\pm k,\pm [k+1])$ has diameter $D=k$.
\end{theorem}

Notice that, for such numbers $N$ of vertices, $k$ is the smallest possible value of the diameter $D$ that solves the $(\Delta,N)$ problem for these graphs. Moreover, for diameter $D=k$, the number $N=2k^2+2k+1$ is the maximum possible number of vertices, so solving the corresponding $(\Delta,D)$ problem.

\subsection{The $(\Delta,N)$ problem for New Amsterdam digraphs}
From the double-step graphs described in Theorem \ref{Bermond-theo}, we solve the  $(\Delta,N)$
problem for the New Amsterdam digraphs for all, but one, (even) values of $N$.

\begin{theorem}
\label{teorema1}
From $a=k$ and $b=k+1$, for $k\ge 1$, we get a New Amsterdam digraph with number $N_{\NA}$ of vertices
and steps
\begin{equation}
\label{steps-NA-theo}
\beta=-\alpha=1,\qquad \gamma=-\delta=2k+1,
\end{equation}
such that
\begin{itemize}
\item[$(a)$]
If $4k^2+4 \le N_{\NA}\le 4k^2+4k+2$, then $D_{\NA}=2k+1$.
\item[$(b)$]
If $N_{\NA}= 4k^2+4k+4$, then $D_{\NA}=2k+2$.
\item[$(c)$]
If $4k^2+4k+8 \le N_{\NA}\le 4(k+1)^2+2$, then $D_{\NA}=2k+3$.
\end{itemize}
\end{theorem}
\prova $(a)$ When $N_{\NA}$ is within the range of $(a)$,  Theorems
\ref{Bermond-theo} and \ref{basic-theo} yield that the diameter is
$D_{\NA}\le 2k+1$ but, from the Moore bounds (\ref{Moore-NA-odd})
and (\ref{Moore-NA-even}), $D_{\NA}$ is smaller than $2k+1$.
Thus, $D_{\NA}=2k+1$. The same reasoning applies for the case
$N_{\NA}= 4(k+1)^2+2$ in $(c)$, giving diameter $D_{\NA}=2k+3$. In
fact if, in this case, we change $k+1$ for $k$ we get the diameter
$D_{\NA}=2k+1$ for the order $N_{\NA}= 4k^2+2$. This is the lower
missing value in $(a)$ to cover all the range $(2k)^2 < N_{\NA}\le
(2k+1)^2+1$ with the minimum possible diameter.

\begin{figure}[h]
\begin{center}
\includegraphics[width=11cm]{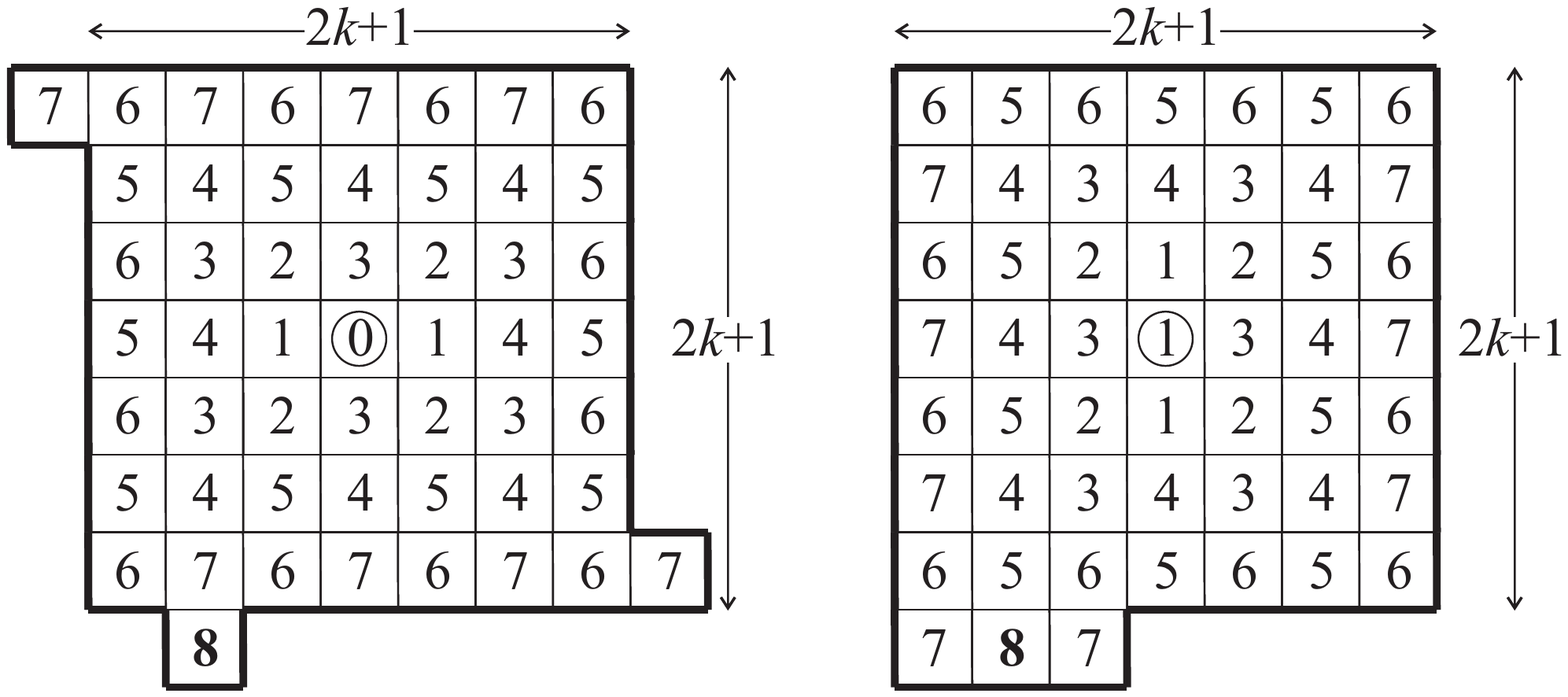}
\vskip-10.75cm
\caption{Optimal tiles for case $(b)$: $N_{\NA}= 4k^2+4k+4$ when $k=3$. In bold, there are the distance of the unique farthest vertex from $0$ or $1$.}
\label{fig-NA-N=52-k=3}
\end{center}
\end{figure}

$(b)$ When $N_{\NA}= 4k^2+4k+4$, we have the tiles centered at $0$
or at $1$ of Fig.~\ref{fig-NA-N=52-k=3},
which gives diameter $D_{\NA}=2k+2$, and their
corresponding tessellations.
Moreover, the distribution of the $0$'s yields the equations
with solutions $a=k$ and $b=k+1$, as claimed (see
Morillo, Fiol, and F\`abrega \cite{mff85}).

$(c)$ Although, according to (\ref{Moore-NA-odd}) and
(\ref{Moore-NA-even}), for the range $4k^2+4k+8 \le N_{\NA}\le
4(k+1)^2$ the minimum theoretical diameter is $D_{\NA}=2k+2$, a
detailed study of the corresponding tiles and tessellations show
that this is not attainable and, then, $D_{\NA}=2k+3$ (a value
guarantied again by Theorems \ref{Bermond-theo} and
\ref{basic-theo}). \final

Note that the only missing value in the above interval
is $N_{\NA}=4k^2+4k+6$. In this case, computer programming
shows that, to obtain the minimum diameter $D_{\NA}=2k+2$, we need to use
other steps different from $a=k$ and $b=k+1$. 


In terms of the diameter $D=D_{\NA}$ of the New Amsterdam digraph, we have the following results:
\begin{itemize}
\item
For every odd diameter $D=D_{\NA}\ge 3$, there is a New Amsterdam digraph for any order $N=N_{\NA}$ satisfying
\begin{equation}\label{resultsNA-odd}
(D-1)^2-2D+10 \le N\le D^2+1,
\end{equation}
to be compared with the theoretical optimal values which, according to the Moore bounds (\ref{Moore-NA-odd}) and (\ref{Moore-NA-even}), would be
\begin{equation}\label{Moore-D-odd-NA}
(D-1)^2+2\le N\le D^2+1.
\end{equation}
\item
For every even diameter $D=D_{\NA}\ge 2$, there is a New Amsterdam digraph for any order $N=N_{\NA}$ satisfying
\begin{equation}\label{resultsNA-even}
D^2-2D+4 \le N\le D^2-2D+6,
\end{equation}
to be compared with the theoretical optimal values which, according to (\ref{Moore-NA-odd}) and (\ref{Moore-NA-even}),  would be
\begin{equation}\label{Moore-D-even-NA}
(D-1)^2+3\le N\le D^2.
\end{equation}
\end{itemize}

Notice that, according to the upper bound in (\ref{resultsNA-odd}), the value $N_{\NA}=4k^2+4k+2$ for odd diameter $D_{\NA}=2k+1$ solves the $(\Delta,D)$ problem (see also Morillo, Fiol and F\`{a}brega~\cite{mff85}).
In fact, for $k=1$ the New Amsterdam digraph $N\!A(10,-1,1,3,3)$
is a bipartite Moore digraph with degree $2$ and diameter $3$. That is, its order attains the Moore bound for a (general) bipartite digraph with such parameters.


\subsection{The $(\Delta,N)$ problem for Manhattan digraphs}
From the above New Amsterdam digraphs, we also can solve the  $(\Delta,N)$
problem for the Manhattan digraphs. The results are based on the already mentioned fact that,
if $G$ is a New Amsterdam digraph on $N_{\NA}$ vertices and diameter $k$, then its line digraph $L(G)$ is a Manhattan digraph on $N_{\MH}=2N_{\NA}$ vertices and diameter $D=k+1$.

\begin{theorem}
\label{teorema2}
From the steps $\beta=-\alpha=1, \gamma=-\delta=2k+1$ of a New Amsterdam digraph, we
get a Manhattan digraph with number $N_{\MH}=2N_{\NA}=4N$ of vertices
and steps
\begin{equation}
\begin{array}{llll}
a_0=1, & a_1=-3, & a_2=1, & a_3=1, \\
b_0=4k+3, & b_1=4k+3, & b_2=-4k-1, & b_3=-4k-5,\label{steps-MH}
\end{array}
\end{equation}
such that
\begin{itemize}
\item
If $8k^2+8 \le N_{\MH}\le 8k^2+8k+4$, then $D_{\MH}=2k+2$.
\item
If $N_{\MH}= 8k^2+8k+8$, then $D_{\MH}=2k+3$.
\item
If $8k^2+8k+16 \le N_{\MH}\le 8(k+1)^2+4$, then $D_{\MH}=2k+4$.
\end{itemize}
\end{theorem}

Now the only missing value of $N_{\MH}$
is $N_{\MH}=8k^2+8k+12$ where, according to the program outputs for the New Amsterdam digraphs, with some other steps different from $\beta=-\alpha=1$ and $\gamma=-\delta=2k+1$, we obtain the diameter $D_{\MH}=2k+3$. 

In terms of the diameter $D=D_{\MH}$ of the Manhattan digraph, we now have the following results:
\begin{itemize}
\item
For every even diameter $D=D_{\MH}\ge 4$, there is a Manhattan digraph for any order $N=N_{\MH}$ satisfying
\begin{equation}\label{resultsMH-even}
2[(D-2)^2-2(D-1)+10] \le N\le 2[(D-1)^2+1],
\end{equation}
to be compared with the theoretical optimal values which, according to the Moore bounds (\ref{Moore-MH-odd}) and (\ref{Moore-MH-even}), would be
\begin{equation}\label{Moore-D-even-NA}
2[(D-2)^2+3]\le N\le 2[(D-1)^2+1].
\end{equation}
\item
For every odd diameter $D=D_{\MH}\ge 5$, there is a Manhattan digraph for any order $N=N_{\MH}$ satisfying
\begin{equation}\label{resultsMH-odd}
2[(D-1)^2-2(D-1)+4] \le N\le 2[(D-1)^2-2(D-1)+6],
\end{equation}
to be compared with the theoretical optimal values which, according to (\ref{Moore-MH-odd}) and (\ref{Moore-MH-even}), would be
\begin{equation}\label{Moore-D-odd-NA}
2[(D-2)^2+3]\le N\le 2(D-1)^2.
\end{equation}
\end{itemize}

Notice that, according to the upper bound in (\ref{resultsMH-even}), the value $N_{\MH}=8k^2+8k+4$ for even diameter $D_{\MH}=2k$ solves the $(\Delta,D)$ problem (see also Morillo, Fiol and F\`{a}brega~\cite{mff85}). In particular, for $k=1$ the Manhattan digraph $M\!H$ with $20$ vertices has diameter $4$, whereas its line digraph $L(M\!H)$, with 40 vertices and diameter 5, turns out to be a quasi-Moore bipartite digraph (that is, it has only two vertices less than the unattainable Moore-bound). For more details about almost Moore bipartite digraphs, see Fiol and Gimbert~\cite{fg01}.

\noindent\textbf{Acknowledgement.}
Research supported by the Ministerio de Educaci\'on y
Ciencia, Spain, and the European Regional Development Fund under
project MTM2011-28800-C02-01 and by the Catalan Research Council
under project 2009SGR1387.


\label{'ubl'}
\end{document}